\documentclass[preprint,12pt]{elsarticle}
\usepackage{amsfonts}
\usepackage{graphics}
\usepackage[
letterpaper, top=1in, bottom=1in, left=1in, right=1in]{geometry}
\usepackage{xcolor}
\usepackage{amssymb}
\usepackage{psfrag}
\usepackage{graphicx}
\usepackage[utf8]{inputenc}
\usepackage{bm}
\usepackage{mathtools}
\usepackage{epstopdf}
\usepackage[colorlinks,linkcolor=black,anchorcolor=black,citecolor=black]{hyperref}
\usepackage{cleveref}
\usepackage{setspace}
\usepackage{amsmath,amssymb,bbm}
\usepackage[utf8]{inputenc}
\usepackage{url}
\usepackage{color}

\newcommand{\be}{\begin{equation}}
\newcommand{\ee}{\end{equation}}
\usepackage{enumitem}

\newtheorem{scheme}{Scheme}

\newproof{proof}{Proof}

\begin{document}
\begin{frontmatter}

\title{An efficient probabilistic scheme for the exit time probability of $\alpha$-stable L\'evy process}

\author[ORNL2]{Minglei Yang}
\author[UTAustin1]{Diego del-Castillo-Negrete}
\author[ORNL1]{Guannan Zhang}
\address[ORNL2]{Fusion Energy Division, Oak Ridge National Laboratory, Oak Ridge, TN.}
\address[UTAustin1]{Institute for Fusion Studies, Dept of Physics, University of Texas at Austin, Austin,  Texas.}
\address[ORNL1]{Computer Science and Mathematics Division, Oak Ridge National Laboratory, Oak Ridge, TN.}

\begin{abstract}
The $\alpha$-stable L\'evy process, commonly used to describe L\'evy flight, is characterized by discontinuous jumps and is widely used to model anomalous transport phenomena. In this study, we investigate the associated exit problem and propose a method to compute the exit time probability, which quantifies the likelihood that a trajectory starting from an initial condition exits a bounded region in phase space within a given time. This estimation plays a key role in understanding anomalous diffusion behavior.
The proposed method approximates the $\alpha$-stable process by combining a Brownian motion with a compound Poisson process. The exit time probability is then modeled using a framework based on partial integro-differential equations (PIDEs). The Feynman-Kac formula provides a probabilistic representation of the solution, involving conditional expectations over stochastic differential equations. These expectations are computed via tailored quadrature rules and interpolation techniques.
The proposed method achieves first-order convergence in time and offers significant computational advantages over standard Monte Carlo and deterministic approaches. In particular, it avoids assembling and solving large dense linear systems, resulting in improved efficiency. We demonstrate the method’s accuracy and performance through two numerical examples, highlighting its applicability to physical transport problems.

\end{abstract}

\tnotetext[fn1]{{\bf Notice}:  This manuscript has been authored by UT-Battelle, LLC, under contract DE-AC05-00OR22725 with the US Department of Energy (DOE). The US government retains and the publisher, by accepting the article for publication, acknowledges that the US government retains a nonexclusive, paid-up, irrevocable, worldwide license to publish or reproduce the published form of this manuscript, or allow others to do so, for US government purposes. DOE will provide public access to these results of federally sponsored research in accordance with the DOE Public Access Plan.}

\end{frontmatter}

\section{Introduction}

Anomalous diffusion plays a critical role in modeling various complex systems across multiple disciplines, including plasma dynamics, mathematical finance, and biological transport~\cite{cartea2007fractional,drummond1962anomalous, hofling2013anomalous}. For instance, it describes nonlocal transport phenomena such as turbulence and heat conduction in magnetized plasmas, where numerical evidence has revealed strongly non-Gaussian particle displacement distributions with algebraic decaying tails that deviate from standard diffusion paradigms \cite{del2004fractional, carreras2001anomalous, del2005nondiffusive}. Unlike normal diffusion, where microscopic time scales are well-defined, anomalous diffusion arises from random walks with divergent microscopic time scales, leading to non-Gaussian stochastic fluctuations \cite{henry2006anomalous, tejedor2010anomalous}. Instead of the linear time dependence seen in normal diffusion, anomalous diffusion is characterized by a power-law scaling of the mean-squared displacement (MSD), $\langle x^2(t) \rangle \sim t^\kappa$, where $\kappa$ indicates the diffusion regime. For $\kappa > 1$, the process is termed superdiffusive, while $0 < \kappa < 1$ corresponds to subdiffusive behavior. The $\alpha$-stable L\'evy process is a prominent example of superdiffusion with scaling $\langle x^2(t) \rangle \sim t^{2/\alpha}$, where $0 < \alpha < 2$ quantifies the degree of nonlocality \cite{meerschaert2006coupled}.

The study of exit problems is fundamental to understanding diffusion-driven transport phenomena. These problems typically focus on determining the first passage or exit time defined as the moment a stochastic process satisfies a given condition or leaves a specified domain for the first time. In the case of normal diffusion, first passage times and related exit statistics are well described by the Fokker-Planck equation with local boundary conditions \cite{condamin2005first, bicout1997green,redner2001guide}.  However, analyzing first passage times for anomalous diffusion is considerably more challenging \cite{gitterman2004reply, chechkin2003first}, and the development of accurate and stable numerical methods for fractional partial differential equations remains an active area of research \cite{lynch2003numerical}. The standard local boundary conditions applicable to Brownian motion no longer hold for the $\alpha$-stable L\'evy process, due to its inherent discontinuous jumps \cite{dybiec2006levy, dybiec2007escape}. Exit problems for symmetric $\alpha$-stable processes with reflecting or absorbing boundaries have been investigated in several works \cite{dybiec2006levy, bass1983exit, deblassie1990first, martin2011first}, and fractional diffusion equations have been applied to plasma transport problems that naturally give rise to such exit time problems \cite{del2004fractional, del2005nondiffusive, del2006fractional}. Moreover, fractional diffusion models have proven effective in describing front dynamics and nonlocal transport phenomena in reaction-diffusion systems with L\'evy flights \cite{del2003front, del2009truncation}.

The primary objective of this paper is to investigate exit problems associated with the L\'evy flight process, in which jump lengths follow a symmetric $\alpha$-stable (S$\alpha$S) distribution characterized by heavy-tailed, power-law decay. Specifically, we focus on computing the probability that particles governed by the $\alpha$-stable process exit a bounded domain before a given time—a quantity referred to as the exit time probability. In the context of random walks, anomalous diffusion, and L\'evy flight behavior, the continuous-time random walk (CTRW) framework has been widely used \cite{scher1973stochastic, klafter1987stochastic}. However, the general CTRW formulation becomes increasingly complicated and less accurate in bounded domains due to the presence of long-range jumps \cite{barkai2000continuous}. As a result, it is not well-suited for studying exit problems involving the $\alpha$-stable L\'evy process. An alternative and commonly used modeling approach involves fractional diffusion equations, which capture the asymptotic behavior of CTRWs and lead to partial integro-differential equations (PIDEs) that provide a comprehensive treatment of anomalous transport phenomena \cite{klafter1987stochastic, ball1987non, klafter1994probability, cartea2007fluid, cont2005finite}. While CTRW-based and fractional PDE approaches are not equivalent in general, their solutions often converge in the asymptotic regime over large spatial and temporal scales. However, the nonlocal nature of fractional derivatives leads to full coefficient matrices in numerical schemes, resulting in $\mathcal{O}(N^2)$ storage and $\mathcal{O}(N^3)$ computational cost, where $N$ is the number of spatial discretization points. These computational challenges become particularly severe when simulating large domains with significant nonlocal interactions \cite{wang2010direct, yuste2005explicit}.

In this study, we propose a novel numerical approach for computing the exit time probability associated with the symmetric $\alpha$-stable process. Our method is based on a Gaussian approximation technique introduced in \cite{asmussen2001approximations, cohen2007gaussian}, which enables us to approximate the $\alpha$-stable process using a stochastic process composed of a Brownian motion and a compound Poisson process. Unlike the original $\alpha$-stable process, both Brownian motion and compound Poisson processes have well-defined analytical probability density functions. In particular, the compound Poisson component features finite jump amplitudes, making the process more tractable for numerical analysis. The exit time probability of the approximated process is governed by a partial integro-differential equation consisting of a local convection-diffusion operator and a nonlocal diffusion operator with an integrable kernel. We impose absorbing boundary conditions, which terminate the particle trajectory upon reaching the domain boundary. To solve this PIDE, we apply the Feynman-Kac formula, which links the terminal-value solution of the PIDE to expectations over forward stochastic differential equations (SDEs). These expectations are evaluated using appropriate quadrature rules and interpolation techniques. In our previous work \cite{yang2021feynman, yang2023probabilistic}, we demonstrated that this method requires significantly fewer floating point operations (FLOPS) compared to standard finite difference schemes—by approximately three orders of magnitude. This efficiency becomes even more distinct in fractional nonlocal problems, where traditional PDE methods must solve large, dense linear systems at every time step. The key contributions of this work are threefold: (i) numerical validation of the Gaussian approximation, showing that the resulting SDE retains the $\alpha$-stable distribution and exhibits the correct power-law scaling of the mean-squared displacement; (ii) the development of a fully discrete numerical scheme tailored for exit problems involving $\alpha$-stable processes; and (iii) numerical examples to verify the accuracy and effectiveness of the proposed numerical scheme for exit problems.

The structure of the paper is as follows. In Section~\ref{sec:problem}, we introduce the definition of the $\alpha$-stable L\'evy process and formulate the exit time probability of interest. Section~\ref{sec:Gauss} presents the Gaussian approximation approach and provides numerical evidence that the resulting process preserves the stable distribution and exhibits the expected power-law behavior in the mean-squared displacement. In Section~\ref{sec:PDE}, we derive the partial integro-differential equation governing the exit time probability and develop a fully discrete probabilistic numerical scheme. Finally, Section~\ref{sec:ex} presents numerical examples to illustrate the effectiveness and accuracy of the proposed method.



\section{Problem setting}\label{sec:problem}

An $\alpha$-stable random variable $\zeta$, with stability index $\alpha \in (0,2]$, is characterized by its characteristic function
\begin{equation}\label{LKeq}
\phi(k) =
\begin{cases}
 \exp\left[ -\sigma^{\alpha} |k|^{\alpha} \left(  1-i \beta \text{sgn}(k)) \tan(\frac{\pi \alpha}{2} \right) + ik\mu \right], &  \alpha \neq 1, \\
 \exp\left[ -\sigma |k| \left(  1 + i \beta  \frac{\pi \alpha}{2} \text{sgn}(k) \ln{|k|}  \right) + ik\mu \right], &\alpha = 1.
\end{cases}
\end{equation}
This defines the $\alpha$-stable distribution $\mathcal{S}(\alpha, \beta, \sigma, \mu)$, where $\alpha$ controls the tail behavior, $\beta \in [-1,1]$ describes the skewness, $\sigma > 0$ is the scale parameter, and $\mu \in \mathbb{R}$ is the location parameter. When $\alpha = 2$, the distribution reduces to the Gaussian; for $\alpha = 1$, it corresponds to the Cauchy distribution.

This work focuses on the $d$-dimensional stochastic differential equation 
\begin{equation}\label{eq_sde}
    X_t^\alpha = X_0^\alpha + \int_{0}^{t} (D_{\alpha})^{1/\alpha} dL_s^{\alpha} \quad \text{with} \quad X_0^\alpha \in \mathcal{D}\subset \mathbb{R}^d,
\end{equation}
where $\alpha \in [1,2]$ (superdiffusion), $\sigma_{\alpha}$ is the diffusion coefficient and  $L_t^{\alpha}$ is a symmetric $\alpha$-stable (S$\alpha$S) L\'evy process. A symmetric $\alpha$-stable L\'evy process is a process with independent increments 
\begin{equation}
    \Delta \mathcal{L}^{\alpha}(\tau)=L^{\alpha}_{t+\tau} - L_t^{\alpha} = \tau^{1/\alpha}  \Delta \mathcal{L}^{\alpha}(1),
\end{equation}
where the increment $\Delta \mathcal{L}^{\alpha}(\tau)$ follows the stable distribution $\mathcal{S}(\alpha,0, \tau^{1/\alpha}, 0)$ and the corresponding characteristic function is defined as 
\begin{equation}
\phi(\tau,w) = \mathbb{E}[e^{iw L_{\tau}^\alpha}] = e^{-\tau|w|^\alpha}.
\end{equation}
According to the definition of the Riesz fractional derivative, the space-fractional diffusion equation associated with the stochastic differential equation (SDE) $X_t^\alpha$ in Eq.~\eqref{eq_sde} can be expressed as 
${\partial \phi}/{\partial t} = D_\alpha {\partial^{\alpha} \phi}/{\partial |x|^{\alpha}}$, where $D_\alpha$ is the diffusion coefficient and the derivative is understood in the sense of the Riesz operator \cite{ibe2013markov, dubkov2008levy}.
The L\'evy  measure describes the frequency and magnitude of jumps in a L\'evy  process. It quantifies the likelihood of jumps of a given size and is a key element in characterizing the discontinuous nature of such processes. For the symmetric $\alpha$-stable (S$\alpha$S) process $L_t^\alpha$, the L\'evy  measure $\nu$ is given by
\begin{equation}\label{levymea}
\nu(dq) = C_{d,\alpha}\frac{1}{|q|^{d+\alpha}}dq, \quad \text{for} \quad q\in \mathbb{R}^{d},
\end{equation}
where the constant 
$
C_{d,\alpha} = \frac{\alpha2^{\alpha-1} \Gamma(\frac{\alpha+d}{2})}{\pi^{d/2}\Gamma(\frac{2-\alpha}{2})}
$ ensures normalization of the jump intensity.
Note that the non-integrable L\'evy measure $\nu(dq)$ satisfies
\begin{equation}
\nu(\{0\}) = 0 \quad \text{and} \quad \int_{\mathbb{R}^d} (1 \wedge|q|^2)\nu(dq) < \infty.
\end{equation}

Let $X_0^\alpha$ be an initial condition (position) of L\'evy process $\{X_t^\alpha, 0\le t\le T_{\rm max}\}$ in an open bounded domain $\mathcal{D} = (\alpha_1, \beta_1) \times \cdots \times (\alpha_d, \beta_d) \subset \mathbb{R}^d$. The first exit time $\tau_X$ of process $X_t^\alpha$ is defined by
\begin{equation}\label{eq_exit_L}
\tau_X := \inf \big\{ t > 0 \,|\, X_0^\alpha \in \mathcal{D}, X_t^\alpha \not\in \mathcal{D}\big\},
\end{equation}
which is the first time instant when the L\'evy process $X_t^\alpha$ exits the domain $\mathcal{D}$.
The goal of this work  is to numerically compute the {\em exit time probability}, $P_X(t,x)$, which is defined as the probability of the event that the stochastic process $X_t^\alpha$ exits the domain $\mathcal{D}$ at or before a specific time  $t \in (0, T_{\max})$, i.e.,
\begin{equation}\label{e2}
P_X(t,x) := \mathbb{P}\left\{ \tau_X \le t \,|\, X_0^\alpha = x \in \mathcal{D} \right\},\; \text{ for any }\; (t,x) \in [0,T_{\max}] \times \overline{\mathcal{D}}.
\end{equation}



Sampling-based methods estimate the exit time probability by simulating trajectories, but obtaining accurate estimates is challenging due to the complex exit behavior of jump processes. L\'evy flights may cross domain boundaries without direct contact, and trajectories may re-enter the domain after exiting. Deterministic alternatives based on fractional Fokker-Planck equations suffer from high computational cost due to nonlocal operators requiring large, dense linear systems.
We propose an alternative strategy based on a Gaussian approximation, constructing an approximate process $\widehat{L}_t^{\alpha}$ for the symmetric $\alpha$-stable process $L_t^\alpha$ using Brownian motion and a compound Poisson process  in Section~\ref{sec:Gauss}. We then develop a fully discrete numerical scheme to compute the exit time probability in Section~\ref{sec:PDE}.


\section{Gaussian approximation of the symmetric $\alpha$-stable process $L_t^\alpha$}\label{sec:Gauss}

The theoretical foundation for the Gaussian approximation has been established in prior works \cite{asmussen2001approximations, cohen2007gaussian}. Early efforts focused on approximating the L\'evy process using an appropriately defined compound Poisson process. However, it is important to note that for $\alpha < 2$, the L\'evy process exhibits paths of infinite variation, which can lead to significant errors if the remainder of the process is neglected. To mitigate this, a Brownian motion with a small variance is introduced to approximate the residual component, thereby improving the fidelity of the approximation.
Specifically, a symmetric $\alpha$-stable (S$\alpha$S) process ${L_t^{\alpha} : 0 \leq t \leq T}$ is approximated by a stochastic process $\widehat{L}_t^{\alpha}$ composed of a compound Poisson process ${N_t^\epsilon : 0 \leq t \leq T}$ with L\'evy measure $\nu^{\epsilon}$ and a Brownian motion ${W_t : 0 \leq t \leq T}$ with small variance $\sigma_\epsilon$, i.e.,
\begin{equation}\label{appro_levy}
L_t^{\alpha} \approx \widehat{L}_t^{\alpha} := ta_{\epsilon} + \sigma_{\epsilon} W_t + N_t^{\epsilon} +  Y_{\epsilon},
\end{equation}
where the $a_{\epsilon}$ is a drift given by
\begin{equation}
a_{\epsilon} = \int_{\|q\|>1} q\nu_{\epsilon}(dz) - \int_{\|q\|\leq1} z\nu^{\epsilon}(dq),
\end{equation}
where $\|\cdot\|$ is the $L^2$ norm. The variance $\sigma_{\epsilon}$ of Brownian motion $W_t$ is defined as $\sigma_{\epsilon} \sigma_{\epsilon}^{\top}=\int_{\mathbb{R}^d} q \cdot q^{\top} \nu_{\epsilon}(dq)$
and the L\'evy measure $\nu$ in Eq.~\eqref{levymea} is decomposed as 
\begin{equation}\label{nu_decom}
\nu=\nu_{\epsilon}+\nu^{\epsilon} :=\nu_{\{\|q\| < \epsilon\}} + \nu_{\{\|q\| \geq \epsilon\}}.
\end{equation}
Under the decomposition in Eq.~\eqref{nu_decom}, it is easy to get $a_{\epsilon}  = 0$ for a symmetric L\'evy process. 
%
%
The kernel $\nu^{\epsilon}(dq)$  of the compound Poisson jump process $N_t^{\epsilon}$ is explicitly expressed as
\begin{equation}\label{eq_kernel}
\nu^{\epsilon}(dq)=\left\{
\begin{array}{rcl}
\begin{aligned}
&C_{d,\alpha}\frac{1}{|q|^{d+\alpha}}dq, \quad {\|q\| \ge \epsilon},\\
&0, \quad \text{otherwise}.
\end{aligned}
\end{array}
\right.
\end{equation}
From Eq.~\eqref{eq_kernel} we have $\nu^{\epsilon}$ is nonnegative and integrable, i.e.,
\begin{equation}\label{cond_lambda}
\nu^{\epsilon}(q) \geq 0 \text{ for } q \in E \quad \text{and}\quad \varphi(q) = \frac{\nu^{\epsilon}(q)} {\lambda} \text{ with } \lambda =\int_{E} \nu^{\epsilon}(q) dq < \infty.
\end{equation}
The function $\varphi(q)$ can be interpreted as a probability density function, with the interaction domain for jumps defined as $E := \{ q \in \mathbb{R}^d, \|q\| \ge \epsilon \}$. Unlike the $S\alpha S$ process $L_t^\alpha$, the compound Poisson process $N_t^{\epsilon}$ features a finite jump intensity $\lambda$, determined by the truncated L\'evy measure over $E$.
The process $Y_\epsilon$, which is c$\grave{a}$dl$\grave{a}$g, accounts for the residual between the exact L\'evy process and its approximation. According to Theorem 3.1 in \cite{cohen2007gaussian}, for any $\epsilon \in (0,1]$, there exists a residual process $Y_\epsilon$ such that the approximation in Eq.~\eqref{appro_levy} holds, and $Y_\epsilon \to 0$ as $\epsilon \to 0$. From a numerical perspective, $\epsilon$ is chosen small enough to ensure the accuracy of the approximation while neglecting $Y_\epsilon$.
However, it is important to note that choosing $\epsilon$ too small leads to a prohibitively high jump intensity $\lambda$ in the compound Poisson process, which can affect numerical stability and efficiency. Based on our numerical experiments, selecting $\epsilon$ in the range $[0.1, 0.3]$ provides a favorable balance. In this regime, the approximating process $\widehat{L}_t^\alpha$ in Eq.~\eqref{eq_xtilde} reliably reproduces the statistical behavior of the target $S\alpha S$ process ${L_t^\alpha : 0 \le t \le T, , \alpha \in [1,2)}$.
Neglecting the residual $Y_\epsilon$, the approximation $\widehat{L}_t^\alpha$ simplifies to the following expression
\begin{equation}\label{eq_ltilde}
\widehat{L}_t^{\alpha} = \sigma_{\epsilon} W_t + N_t^{\epsilon}.
\end{equation}

We now clarify the distinction between the trajectories of the $S\alpha S$ process $L_t^{\alpha}$ and its Gaussian approximation $\widehat{L}_t^{\alpha}$, as defined in Eq.~\eqref{eq_ltilde}. The path of $L_t^{\alpha}$ is almost surely discontinuous at all times and typically exits a bounded domain $\mathcal{D}$ without directly interacting with its boundary. In contrast, the trajectory of $\widehat{L}_t^{\alpha}$ exhibits discontinuities only at discrete jump times and may exit the domain either by directly hitting the boundary or by landing in the exterior $D^c$ due to a jump.
In this work, our focus lies not on reproducing microscopic trajectory behavior but rather on ensuring that both processes, $L_t^{\alpha}$ and $\widehat{L}_t^{\alpha}$, exhibit consistent statistical properties—such as the mean-squared displacement and probability density function. This level of agreement is sufficient for accurately approximating the exit time probability and serves our numerical objective well.

\subsection{Numerical validation of Gaussian approximation}\label{sec:num_gauss}
This section aims to numerically verify the accuracy of the Gaussian approximation. Specifically, we demonstrate the agreements on two statistical perspectives:
\begin{itemize}
    \item Both processes $L_t^{\alpha}$ and $\widehat{L}_t^{\alpha}$ follow the $\alpha$-stable distribution $\mathcal{S}(\alpha, 0, t^{1/\alpha}, 0)$.
    \item Both processes $L_t^{\alpha}$ and $\widehat{L}_t^{\alpha}$ exhibit a growth rate of $2/\alpha$ in mean-squared displacement (MSD), denoted as $\langle x^2(t) \rangle \sim t^{2/\alpha}$.
\end{itemize}

\begin{figure}[h!]
    \centering
   {\includegraphics[scale = 0.5]{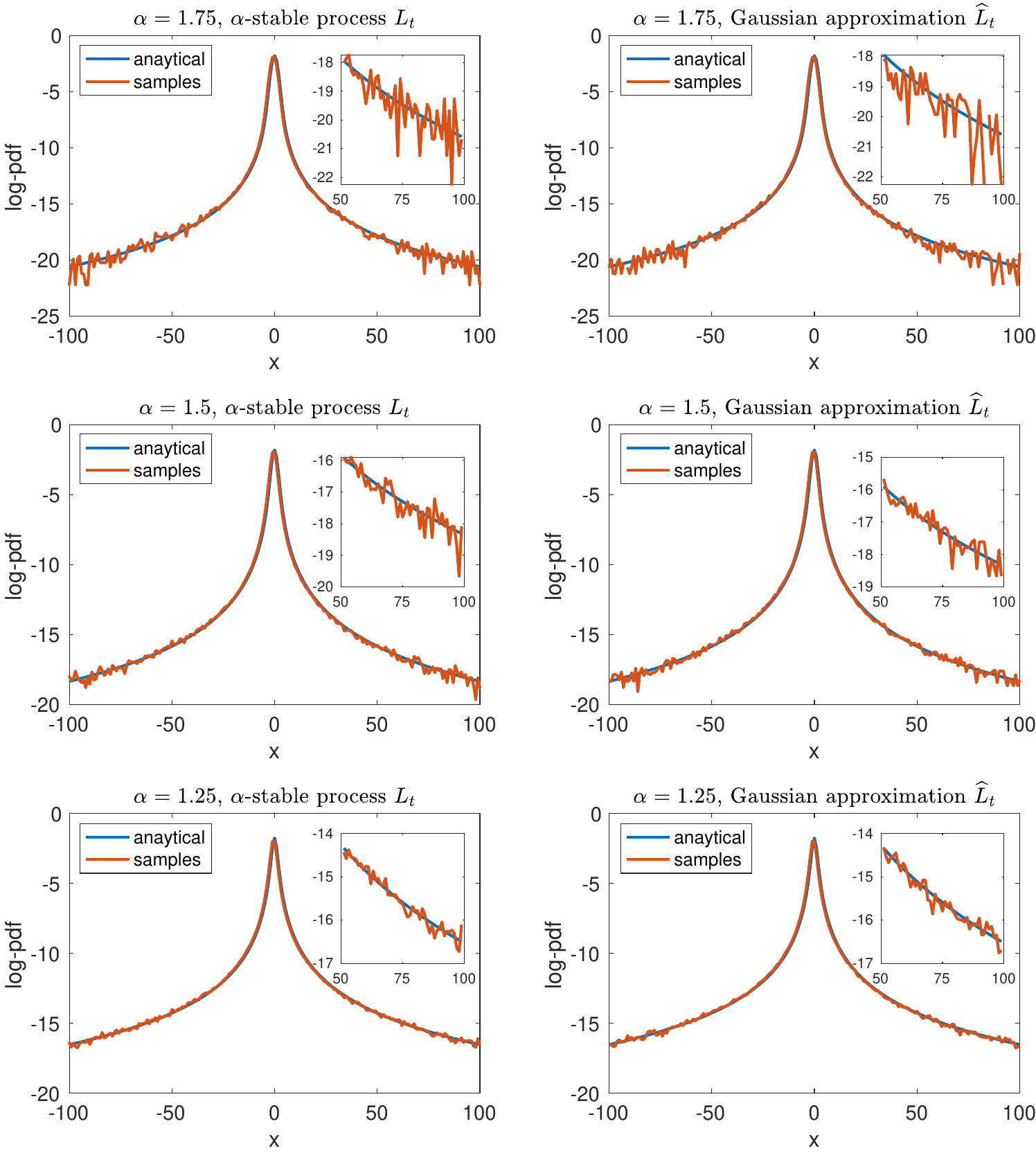} }
    \caption{The logarithmic probability density functions (log-pdf) of two processes, $L_t^{\alpha}$ and its approximation $\widehat{L}_t^{\alpha}$, at the unit time $t = 1$. The study tests three distinct $\alpha$ values: $1.75$, $1.5$, and $1.25$, each presented in a separate row. The blue curve represents the theoretical log-pdf produced through the MATLAB STBL toolbox \cite{stable}, while the red curve depicts the log-pdf based on $N = 10^5$ samples. Both processes, $L_t^{\alpha}$ and $\widehat{L}_t^{\alpha}$, exhibit accurate adherence to the $\alpha$-stable distribution, as highlighted in the figure.}
    \label{dist_approx}
\end{figure}

Fig.~\ref{dist_approx} presents the logarithmic probability density functions (log-pdfs) of the $L_t^{\alpha}$ and $\widehat{L}_t^{\alpha}$ processes at time $t = 1$, for three representative values of $\alpha$: $1.75$, $1.5$, and $1.25$. Each row of the figure corresponds to a specific value of $\alpha$. The analytical log-pdf, computed using the MATLAB STBL toolbox \cite{stable}, is shown as a blue curve. The red curves represent the empirical log-pdfs obtained from $N = 10{,}000$ samples of both the $\alpha$-stable process $L_t^{\alpha}$ (displayed in the left column) and its Gaussian approximation $\widehat{L}_t^{\alpha}$ (displayed in the right column).
As illustrated in Fig.~\ref{dist_approx}, both $L_t^{\alpha}$ and $\widehat{L}_t^{\alpha}$ closely match the theoretical $\alpha$-stable distribution, confirming the validity of the Gaussian approximation at the level of marginal distributions.

\begin{figure}[h!]
    \centering
  {\includegraphics[scale = 0.35]{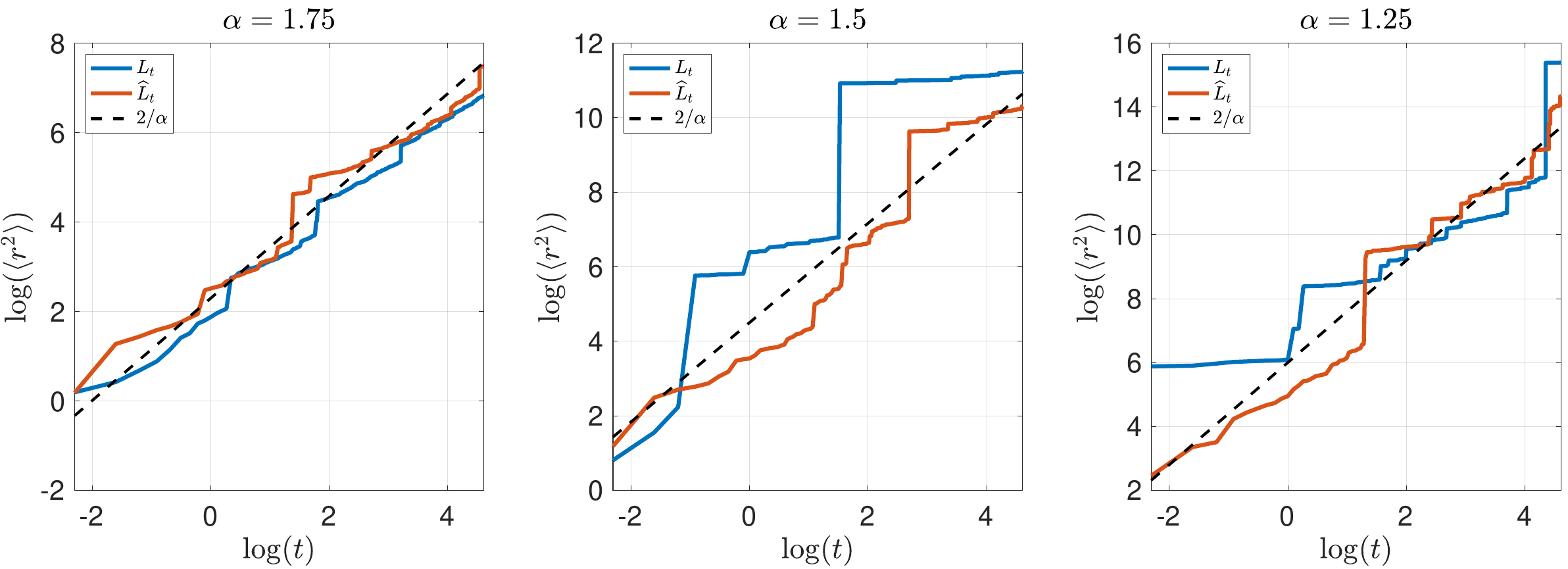} }
    \caption{This figure presents the mean-squared displacement (MSD) plots for the $S\alpha S$ process $L_t^{\alpha}$ and its Gaussian approximation $\widehat{L}_t^{\alpha}$, with separate panels corresponding to $\alpha = 1.75$, $1.5$, and $1.25$. The horizontal axis shows $\log(t)$, and the vertical axis shows $\log(\langle r^2 \rangle)$. Both processes exhibit the expected superdiffusive scaling behavior, with MSD growing as $t^{2/\alpha}$. The approximation $\widehat{L}_t^{\alpha}$ displays reduced fluctuations due to jump truncation but retains the same scaling rate, confirming that it preserves the key statistical feature of anomalous diffusion.}
    \label{fig:MSD}
\end{figure}

Fig.~\ref{fig:MSD} presents the mean-squared displacement (MSD) plots for the processes $L_t^{\alpha}$ and $\widehat{L}_t^{\alpha}$. The panels on the left, middle, and right correspond to $\alpha = 1.75$, $1.5$, and $1.25$ respectively. The $x$-axis represents the logarithm of time ($\log(t)$), while the $y$-axis represents the logarithm of MSD ($\log(\langle r^2 \rangle)$). As depicted in Fig.~\ref{fig:MSD}, both $L_t^{\alpha}$ and $\widehat{L}_t^{\alpha}$ processes exhibit a consistent growth rate of $2/\alpha$ in MSD. Notably, the MSD of $\widehat{L}_t^\alpha$ displays reduced fluctuations due to the truncation of the jump amplitude in the compound Poisson process, where the upper bound of jump amplitude is truncated by $\epsilon_{\rm out} = 10^5$. Despite this, it maintains the same $2/\alpha$ growth rate in MSD. This collectively signifies the superdiffusive nature of the $\widehat{L}_t^{\alpha}$ process.
We set $\epsilon = 0.1$, and use a temporal step size of $\Delta t = 0.005$ in the MSD simulations. 
The combined insights from Figs.~\ref{dist_approx} and \ref{fig:MSD} provide strong evidence that the Gaussian approximation process $\widehat{L}_t^{\alpha}$ faithfully reproduces key statistical properties of the original $S\alpha S$ process $L_t^{\alpha}$. Specifically, it follows the $\alpha$-stable distribution and preserves the expected mean-squared displacement (MSD) growth rate of $t^{2/\alpha}$. 

Although the second moment of an $\alpha$-stable L\'evy process with $\alpha < 2$ is formally infinite, the mean-squared displacement (MSD) shown here should be interpreted as an effective (truncated) second moment of the approximated process, and is used solely to verify the expected superdiffusive scaling exponent $2/\alpha$.

\section{Numerical scheme for the exit time probability of $\widehat{X}_t^{\alpha}$}\label{sec:PDE}

Based on the Gaussian approximation introduced in Section.~\ref{sec:Gauss}, we obtain the approximation of the SDE of interest $X_t^\alpha$ in Eq.~\eqref{eq_sde} as follows
\begin{equation}\label{eq_xtilde}
{X}_t^\alpha  \approx \widehat{X}_t^\alpha = \widehat{X}_0^\alpha + \int_{0}^t (D_\alpha)^{1/\alpha} \sigma_{\epsilon}dW_s +  (D_\alpha)^{1/\alpha} N_t^\epsilon\quad \text{with} \quad \widehat{X}_0^\alpha = X_0^\alpha. 
\end{equation}

We now define the exit time probability associated with the Gaussian approximation $\widehat{X}_t^\alpha$. 
Similar to the definition for the original $S\alpha S$ process $X_t^\alpha$, the exit time probability of the process $\widehat{X}_t^\alpha$, denoted by $P_{\widehat{X}^\alpha}(t, x)$, is defined as
\begin{equation}\label{e2_2}
P_{\widehat{X}^\alpha}(t,x) = \mathbb{P}\left\{ \tau_{\widehat{X}} \le t \,|\, \widehat{X}_0^{\alpha} = x \in \mathcal{D} \right\}\; \text{ for any }\; (t,x) \in [0,T_{\max}] \times \overline{\mathcal{D}},
\end{equation}
where the first exit time $\tau_{\widehat{X}^\alpha}$ of $\widehat{X}_t^{\alpha}$ is defined as
\begin{equation}\label{eq_exit_X}
\tau_{\widehat{X}^\alpha} := \inf \big\{ t > 0 \,|\, \widehat{X}_0^{\alpha} \in \mathcal{D}, \widehat{X}_t^{\alpha} \not\in \mathcal{D}\big\}.
\end{equation}

By applying the It\^o formula for jump processes \cite{applebaum2009levy}, extending Theorem 1.4.5 from \cite{Schuss:2013th}, we establish that for any fixed time $T \in [0, T_{\max}]$, the exit time probability $P_{\widehat{X}^\alpha}(T, x)$ defined in Eq.~\eqref{e2_2} can be obtained from the solution of the following equation
\begin{equation}\label{e5}
P_{\widehat{X}^\alpha}(T, x) = u(0, x).
\end{equation}
The solution $u(t, x)$ satisfies the terminal boundary value problem,
\begin{equation}\label{e3}
\begin{aligned}
\frac{\partial u(t,x)}{\partial t}  + \mathcal{L}[u](t,x) & = 0 \quad \text{ for }\;\; x \in \mathcal{D},\, t < T,\\
u(t,x) & = 1  \quad \text{ for }\;\; x \in \mathcal{D}^c,\, t<T,\\[4pt]
 u(T,x) & = 0  \quad \text{ for }\;\; x \in \mathcal{D},\\
\end{aligned}
\end{equation}
where the operator $\mathcal{L}(t,x)$ is
the backward Kolmogorov operator associated with $\widehat{X}_t^{\alpha}$ in Eq.~\eqref{eq_xtilde} given by
\begin{equation}\label{eq:operatorL}
    \mathcal{L}[u](t,x) := \frac{1}{2}\sum_{i,j=1}^{d}(\sigma\sigma^{\top})_{i,j}\frac{\partial^2 u}{\partial x^i \partial x^j}(t,x) +\int_{E}\bigg(u(t,x+c(q))-u(t,x)\bigg)\nu^{\epsilon}(dq),
\end{equation}
where $\sigma = (D_\alpha)^{1/\alpha} \sigma_{\epsilon}$ is the local diffusion coefficient
with $(\sigma\sigma^{\top})_{i,j}$ the $(i,j)$-th entry of the diffusivity tensor $\sigma\sigma^{\top}$, and $c(q) =  (D_\alpha)^{1/\alpha} q$ is the jump amplitude.
In this case, the nonlocal integral component in $\mathcal{L}$ corresponds to the compound Poisson process $N^{\epsilon}_t$, where the interaction domain $E$ and jump kernel $\nu^{\epsilon}$ are defined in Eqs.~\eqref{eq_kernel}. 
The well-posdeness of the problem in Eq.~\eqref{e3} has been discussed and proved in \cite{du2012analysis,du2014nonlocal}.

In the remainder of Section \ref{sec:PDE}, we formulate a probabilistic numerical scheme to solve the terminal boundary value problem in Eq.~\eqref{e3}. Our approach leverages the Feynman-Kac formula to express the solution $u(t,x)$ as a conditional expectation. Using a temporal-spatial mesh, the numerical method achieves first-order convergence in time. We focus our convergence analysis on the exit time probability $P_{\widehat{X}^\alpha}$ for the approximated process rather than the original process $P_{{X}^\alpha}$, as the convergence properties of the underlying Gaussian approximation have been thoroughly analyzed in existing literature \cite{asmussen2001approximations, cohen2007gaussian}.

\subsection{The Feynman-Kac representation}\label{Feynman}
We introduce a uniform mesh over the temporal domain $[0,T]$ as follows
\begin{equation}\label{t_mesh}
\mathcal{T} := \{0=t_{0} \le \cdots \le t_{N_t}=T\},
\end{equation}
with $\Delta t= t_{n}-t_{n-1}$, for $1\leq n \leq N_t$. In each sub-interval $[t_n,t_{n+1}]$, we rewrite $\widehat{X}_t^{\alpha}$ in Eq.~\eqref{eq_xtilde} into the condition form
\begin{equation}\label{sde} 
\widehat{X}_{s}^{t_n} = x + \sigma \Delta W_{s-t_n} + \sum_{k = 1}^{N_{s-t_{n}}} c(q_k),
\end{equation}
where diffusion and jump coefficients $\sigma$ and $c(q)$ are given in Eq.~\eqref{eq:operatorL}, the superscript ${}^{t_n}$ indicates $\widehat{X}_s^{t_n}$ starts from the location $(t_n,x)$ and move forward to $s\in[t_n,t_{n+1}]$ and the Brownian motion increment $\Delta W_{s-t_n} := W_{s} - W_{t_n}$ satisfies $\mathbb{V}ar(W_t) = dt$. $N_{s-t_{n}}$ is the Poisson process following the Poisson probability distribution $
\mathbb{P}(N_{s-t_n}=k) = (\lambda (s-t_n))^k {e^{-\lambda (s-t_n)}}/{k!},
$
with $\lambda$ defined in Eq.~\eqref{cond_lambda}, $t_k$ for $k = 1, \ldots N_{s-t_{n}}$ is the instance of time that jumps occur, and $q_k$ follows the probability distribution defined by $\varphi(q)$ in Eq.~\eqref{cond_lambda}. 
Accordingly, we can define the first conditional exit time 
\begin{equation}
    \tau_n:= \inf \big\{ s > t \,|\, x \in \mathcal{D}, \widehat{X}_s^{t_n} \not\in \mathcal{D}\big\}.
\end{equation}
The probabilistic representation of $u(t_n,x)$ in Eq.~\eqref{e3} is then given by the Feynman-Kac formula as 
\begin{equation}\label{prob_rep}
    u(t_{n},x) = \mathbb{E}\left[ u(t_{n+1}\wedge \tau_{n}, \widehat{X}_{t_{n+1}\wedge \tau_{n}}^{t_n}) \right],
\end{equation}
where $t_{n+1}\wedge \tau_{n}:= \min (t_{n+1},\tau_n)$ denotes the minimum of $t_{n+1}$ and $\tau_n$.
Instead of solving the adjoint equation in Eq.~\eqref{e3} directly, our approach is to  approximate the solution $u$ at each time step by discretizing the Feynman-Kac representation in Eq.~\eqref{prob_rep}. The details of this approximation will be described in Section \ref{sec:approxfey}.

\subsection{The approximation of the Feynman-Kac representation}\label{sec:approxfey}
This section aims to approximate the representation of $u(t_n,x)$ in Eq.~\eqref{prob_rep}. The approximation mainly consists of three steps: (a) numerical treatment of the exit time $\tau_n$; (b) approximation of the expectation $\mathbb{E}[\cdot]$ in Eq.~\eqref{prob_rep}; (c) reconstruction of $u(t_n,x)$ in the spatial domain $\mathcal{D}$. These three steps will be discussed in Section \ref{sec:exit}, \ref{sec:exp}, and \ref{sec:interp}. And a fully discrete scheme is presented in Section \ref{sec:full}.

\subsubsection{Treatment of the exit time}\label{sec:exit}
This section explains how we handle the exit time $\tau_n$ in Eq.~\eqref{prob_rep}. Dealing with exit times becomes challenging when approximating mathematical expectations, especially with jump processes involved, due to the complex exit behaviors. Common methods, like truncating $\tau_n$ to $t_n$, typically result in a half-order convergence rate concerning $\Delta t$ \cite{Gobet:2000cj, buchmann2003computing}. However, in this scenario, we've developed a straightforward approach tailored for nonlocal problems. This treatment achieves a reliable first-order convergence rate with respect to $\Delta t$.
To proceed, the expectation $\mathbb{E}\left[ u( t_{n+1} \wedge \tau_n, \widehat{X}^{t_n}_{t_{n+1} \wedge \tau_n}) \right]$ in Eq.~\eqref{prob_rep} can be decomposed based on different scenarios of $N_{\Delta t} = N_{t_{n+1} - t_n}$ as follows 
\begin{equation}\label{eq_expand}
\begin{aligned}
&  \mathbb{E}\left[ u(t_{n+1}\wedge \tau_{n}, \widehat{X}_{t_{n+1}\wedge \tau_{n}}^{t_n}) \right]= I_1 + I_2 + I_3,\\[2pt]
& I_1 = \mathbb{P}(N_{\Delta t} = 0) \mathbb{E}\left[ u(t_{n+1}\wedge \tau_{n}, \widehat{X}_{t_{n+1}\wedge \tau_{n}}^{t_n}) | N_{\Delta t} = 0\right],\\[2pt]
& I_2 = \mathbb{P}(N_{\Delta t} = 1) \mathbb{E}\left[ u(t_{n+1}\wedge \tau_{n}, \widehat{X}_{t_{n+1}\wedge \tau_{n}}^{t_n}) |N_{\Delta t} = 1\right],\\[2pt]
& I_3 = \sum_{k=2}^{\infty}\mathbb{P}(N_{\Delta t}  = k) \mathbb{E}\left[ u(t_{n+1}\wedge \tau_{n}, \widehat{X}_{t_{n+1}\wedge \tau_{n}}^{t_n})|N_{\Delta t} =k \right],
\end{aligned}
\end{equation}
where each term is an expectation conditional on an event defined by $N_{\Delta t}$.

Moving forward, we analyze the elements in Eq.~\eqref{eq_expand} to identify those that can be treated as negligible in the numerical process. Among the terms we choose to retain in the final numerical approach, we aim to avoid directly approximating the exit time. Specific treatment of each term is given as follows:

\begin{itemize}[leftmargin=20pt]\itemsep0.15cm
\item When $N_{\Delta t} = 0$, the SDE $\widehat{X}_t^{t_n}$ in Eq.~\eqref{sde} is only driven by Brownian motion. It is known that $\mathbb{P}(N_{\Delta t} = 0,\tau_{n} \le t_{n+1}) \rightarrow 1$ as $x\rightarrow \partial \mathcal{D}$. In our previous work \cite{yang2021feynman}, we proved that if the starting location $x$ of $\widehat{X}_s^{t_n}$ satisfies  
\begin{equation}\label{cond_bound}
{\rm dist}(x, \partial \mathcal{D}) \ge \mathcal{O}((\Delta t)^{\frac{1}{2}-\varepsilon}),
\end{equation}
for an arbitrarily small positive number $\varepsilon>0$ with ${\rm dist}(\cdot, \cdot)$ denoting the Euclidean distance, then for sufficiently small $\Delta t$, 
\begin{equation}\label{stoperr}
\mathbb{P}(N_{\Delta t} = 0,\tau_{n} \le t_{n+1}) \leq C (\Delta t)^\varepsilon \exp\left(-\frac{1}{(\Delta t)^{2\varepsilon}}\right),
\end{equation}
where the positive constant $C$ is independent of $\Delta t$. The condition in Eq.~\eqref{cond_bound} can be satisfied by properly defining the spatial mesh, which will be discussed in Section \ref{sec:interp}. 
The estimate in Eq.~\eqref{stoperr} allows us to estimate $I_1$ as
\begin{equation}
    I_1 = \mathbb{P}(N_{\Delta t} = 0) \mathbb{E}\left[ u(t_{n+1}, \widehat{X}_{t_{n+1}}^{t_n}) | N_{\Delta t} = 0\right] + \mathcal{O}((\Delta t)^2).
\end{equation}

\item In terms of $I_2$, the SDE $\widehat{X}_t^{t_n}$ has one jump activity in each time interval $[t_n,t_{n+1})$ along with Brownian motion. As $\Delta t \rightarrow 0$, the orbit's increment from the jump process, which is independent with $\Delta t$, will dominate the increment from Brownian motion.
We introduce the auxiliary variable $V_x^{t_n}$, defined as 
\begin{equation}
    V_s^{t_n} := \widehat{X}_{s}^{t_n} - \sigma \Delta W_{s-t_n} = x + \sum_{k = 1}^{N_{s-t_{n}}} c(q_k),
\end{equation}
which is the truncation of variable $\widehat{X}_s^{t_n}$ increment by only keeping the jump component. Lemma 4.4 in our previous work \cite{yang2023probabilistic} proved that when the condition in Eq.~\eqref{cond_bound} is satisfied, we can 
replace process $\widehat{X}_s^{t_n}$ by $V_s^{t_n}$ and eliminate the conditional exit time $\tau_n$ as follows
\begin{equation}\label{I2_approx}
    I_2 = \mathbb{P}(N_{\Delta t} = 1) \mathbb{E} \left[  u(t_{n},{V}_{t_{n+1}}^{n})\, \big|\,N_{\Delta t} = 1\right] + \mathcal{O}((\Delta t)^2).
\end{equation}
The local truncation error, denoted as $\mathcal{O}((\Delta t)^2)$, in Eq.~\eqref{I2_approx} comprises two components. First, the probability of having a single jump ($N_{\Delta t} = 1$) is of the order $\mathcal{O}(\Delta t)$. Second, the error arising from the elimination of the exit time $\tau_n$, due to martingale property of Brownian motion, is also on the order of $\mathcal{O}(\Delta t)$. Altogether, we achieve the desired local error of $\mathcal{O}((\Delta t)^2)$.

\item Finally, for $I_3$, the probability of the Poisson process $N_{\Delta t}$ having $k$ jumps within $[t_n,t_{n+1})$ is on the order of $\mathcal{O}((\Delta t)^k)$, we have 
$
I_3 = \mathcal{O}((\Delta t)^2)
$
when $\mathbb{E} [u(t_{n+1}\wedge \tau_n, \widehat{X}^{n}_{t_{n+1}\wedge \tau_n})|N_{\Delta t} =k]$ for $k\ge 2$. So we can neglect $I_3$ in the final numerical scheme.

\end{itemize}

Based on the above estimates, we approximate $u(t_{n},x)$ in Eq.~\eqref{prob_rep} by neglecting all truncation errors those are on the order of $\mathcal{O}((\Delta t)^2)$ as 
\begin{equation}\label{qpprox_exit}
\begin{aligned}
 u(t_{n},x) \approx u^n(x) =&\; \mathbb{P}(N_{\Delta t} = 0) \mathbb{E}\left[u^{n+1}(\widehat{X}_{t_{n+1}}^{t_n})\, \big|\, N_{\Delta t} = 0\right] \\
 &\;+ \mathbb{P}(N_{\Delta t} = 1) \mathbb{E} \left[  u^{n+1}({V}^{t_n}_{t_{n+1}})\, \big|\,N_{\Delta t} = 1\right].
\end{aligned}
\end{equation}
with
\begin{equation}\label{Euler}
\widehat{X}_{t_{n+1}}^{t_n} = x + \sigma \,\Delta W, \quad {V}^{t_n}_{t_{n+1}} = x + c(q_1),
\end{equation}
where $\Delta W = W_{t_{n+1}} - W_{t_{n}}$, and jump amplitude $q_1$ follows the probability distribution $\varphi(q)$ defined in Eq.~\eqref{cond_lambda}.

\subsubsection{Approximation of the conditional expectations}\label{sec:exp} 
The computation involving expectation $\mathbb{E}[\cdot]$ in Eq.~\eqref{qpprox_exit} is defined in the entire domain $\mathcal{D}$, making their direct calculation quite challenging.
To overcome this, we establish adequate quadrature rules to approximate the conditional expectations $\mathbb{E}[\cdot]$ in Eq.~\eqref{qpprox_exit}. The expectation $\mathbb{E}[u^{n+1}(\widehat{X}_{t_{n+1}}^{t_n})\, \big|\, N_{\Delta t} = 0]$ only involves the Brownian motion due to $N_{\Delta t} = 0$, i.e.,
\begin{equation}\label{E0}
\mathbb{E}\left[u^{n+1}(\widehat{X}_{t_{n+1}}^{t_n})\, \big|\, N_{\Delta t} = 0\right]=
\int_{\mathbb{R}^d} u^{n+1}(x + \sigma\sqrt{2\Delta t}\xi)\rho(\xi)d\xi,
\end{equation}
where $\xi := (\xi_1, \ldots, \xi_d)$ follows the normal distribution with the probability density 
$\rho(\xi) := \pi^{-d/2}  \exp(-\sum_{\ell=1}^d \xi_\ell^2)$.
We use the tensor-product Gauss-Hermite quadrature rule to approximate this integral, and denote the approximate expectation as 
\begin{equation}\label{E1}
    \widehat{\mathbb{E}}\left[u^{n+1}(\widehat{X}_{t_{n+1}}^{t_n})\, \big|\, N_{\Delta t} = 0\right] := \sum_{m=1}^{M}w_{m}\,u^{n+1}\big(x + \sigma\sqrt{2\Delta t}\, e_{m}\big),
\end{equation}
where $\{w_{m}, e_{m}\}_{m=1}^M$ denote the Gauss-Hermite quadrature weights and abscissa\footnote{We use a single index to represent the tensor-product quadrature rule.}. 

The expectation $\mathbb{E} [u^{n+1}({V}^{t_n}_{t_{n+1}}) \; \big|\; N_{\Delta t} = 1]$ only involves the compound Poisson jumps, i.e.,
\begin{equation}\label{int_ejump}
\mathbb{E}\left[u^{n+1}({V}_{t_{n+1}}^{t_n})\, \big|\, N_{\Delta t} = 1\right]=
\int_{E} u^{n+1}(x + c(q))\varphi(q) dq,
\end{equation}
where $\varphi(q)$ is defined in Eq.~\eqref{cond_lambda}. In this case, the choice of the quadrature rule is determined by $\varphi(q)$. For example, if $\varphi(q)$ is bounded and has a compact support, we can use a Gauss-Legendre rule or a  Newton-Cotes rule; if $\varphi(q)$ is singular at the origin, e.g., $\varphi(q) = 1/|q|^z$ with $0<z<1$, then we can use a Gauss-Jacobi rule. 
In general, we write the quadrature approximation of the compound Poisson jump as
\begin{equation}\label{E2}
   \widetilde{\mathbb{E}}\left[u^{n+1}({V}_{t_{n+1}}^{t_n})\, \big|\, N_{\Delta t} = 1\right] := \sum_{l=1}^{L} v_{l}\, u^{n+1}\big(x + c(a_l)\big),
\end{equation}
where $\{v_{l}, a_l\}_{l = 1}^L$ denote the corresponding quadrature weights and abscissa. 

We now update the approximation $u^n(x)$ by substituting the quadrature rules into Eq.~\eqref{qpprox_exit}, i.e.,
\begin{equation}\label{eq_utn_approx_tilde}
 u^n(x) = \mathbb{P}(N_{\Delta t} = 0) \widehat{\mathbb{E}}\left[u^{n+1}(\widehat{X}_{t_{n+1}}^{t_n})\, \big|\, N_{\Delta t} = 0\right]  + \mathbb{P}(N_{\Delta t} = 1) \widetilde{\mathbb{E}} \left[  u^{n+1}({V}^{t_n}_{t_{n+1}})\, \big|\,N_{\Delta t} = 1\right].
\end{equation}
The accuracy of quadrature rules can be controlled within an error of $\mathcal{O}((\Delta t)^2)$ by using an adequate number of quadrature points. One case involves using the Gauss-Hermite quadrature rule for the Brownian motion, while the trapezoidal rule is employed to approximate the integral associated with the jump component.

Let $M$ denote the number of Gauss-Hermite quadrature points in each dimension. If $u(t,\cdot)$ is sufficiently smooth, i.e., $\partial^{2M} u/\partial \xi^{2M}$ is bounded, then the Hermite quadrature error is bounded by \cite{2013JSV...332.4403B,yang2021feynman}
\begin{equation}\label{GH_quad}
\left|{\mathbb{E}}\left[u^{n+1}(\widehat{X}^{t_n}_{t_{n+1}})\,|\,N_{\Delta t} = 0\right] - \widehat{\mathbb{E}}\left[u^{n+1}(\widehat{X}^{t_n}_{t_{n+1}})\,|\,N_{\Delta t} = 0\right]\right| \le C\frac{M!}{2^{M}(2M)!} (\Delta t)^{M},
\end{equation}
where the constant $C$ is independent of $M$ and $\Delta t$. The factor $(\Delta t)^{M}$ comes from the $2{M}$-th order differentiation 
of $u(t,\cdot)$ with respect to $\xi$ defined in Eq.~\eqref{E0}.
To approximate the integral with respect to jump variable $q$ in Eq.~\eqref{int_ejump}, we divide the interaction domain $E$ by equally spaced mesh size $h>0$. 
Using trapezoidal rule in Eq.~\eqref{E2}, we have the bound
\begin{equation}\label{NC_quad}
\left|\mathbb{E}[u^{n+1}({V}^{t_n}_{t_{n+1}})\, \big|\, N_{\Delta t} = 1] - \widetilde{\mathbb{E}}[u^{n+1}({V}^{t_n}_{t_{n+1}})\, \big|\, N_{\Delta t} = 1]\right|\le C h^2,
\end{equation}
where constant $C$ depends on the second derivative $\partial^{2} u/\partial q^{2}$.
Combining Eq.~\eqref{GH_quad} and Eq.~\eqref{NC_quad}, with $M\ge 2$ and $h \le \sqrt{\Delta t}$, the truncation error from the quadrature rules can be bounded by $\mathcal{O}((\Delta t)^2)$.

In the numerical scheme, we do not need to define the quadrature rule in the whole large-volume jump interaction domain $E$ (which is an unbounded domain without truncation). Beyond the interaction domain $\mathcal{D}$, the exit time probability equals one, i.e., $u(t,x) = 1$ for $x\in \mathcal{D}^c$. Consequently, the quadrature rule is only needed  within the domain $\mathcal{D}$. This outcome has important implications for the implementation of the quadrature method for the exit problems within the bounded domain.

\subsubsection{Spatial approximation}\label{sec:interp}
To approximate $u^n(x)$ within the domain $\mathcal{D}$, we employ piecewise polynomial interpolation. Initially, we introduce a Cartesian mesh $\mathcal{S} = \mathcal{S}^1 \times \mathcal{S}^2 \times \cdots \times \mathcal{S}^d$ for the extended domain $\overline{\mathcal{D}}$, where $\mathcal{S}^i$ for $i = 1, \ldots, d$, represents a mesh defined on the interval $[\alpha_i, \beta_i]$. It is worth noting that the partition along each dimension can be non-uniform.
For notational simplicity, we use $J$ to denote the total number of grid points in $\mathcal{S}$, and use a scalar $j$ to index all the grid points in $\mathcal{S}$, i.e.,
$
\mathcal{S} := \{ x_j: \; j = 1, \ldots, J\}.$
Note that we use a single index to denote the grid points $x_j$ to simplify the notation. 
In particular, we define the approximation of $u^n( x)$ using a $p$th order Lagrange nodal basis as
\begin{equation}\label{eq_interp_E}
  u^{n,p}(x) := \sum_{j =1}^J u^n(x_j)\psi_j(x),
\end{equation}
where $\psi_j$ is the nodal basis function associated with the grid point $x_j$, and $u^n(x_j)$ is the nodal value at $x_j$. 
Since $u(t, x)$ represents a probability value bounded within $[0,1]$, we adopt the Piecewise Cubic Hermite Interpolating Polynomial (PCHIP) method \cite{fritsch1980monotone}, which ensures monotonicity and avoids overshoots, making it particularly suitable for probabilistic functions.

Recall that the estimate in Eq.~\eqref{stoperr} requires the condition in Eq.~\eqref{cond_bound} imposed on the starting location of ${X}_{s}^{t_n}$. This condition is realized by letting the spatial mesh $\mathcal{S}$ satisfy  
\begin{equation}\label{cond_bound1}
{\rm dist}(x_j, \partial \mathcal{D}) \ge \mathcal{O}((\Delta t)^{\frac{1}{2}-\epsilon}) \;\text{ for }\; x_j \in \mathcal{S} \cap \mathcal{D}.
\end{equation}
In actuality, it suffices to apply this condition only to the layer of grid points in close to the boundary $\partial \mathcal{D}$.
In the practical implementation, we achieve this by constructing the mesh $\mathcal{S}$ in such a way that the quadrature points employed in Eq.~\eqref{E1} for all interior grid points invariably reside within the domain. $\mathcal{D}$, i.e., 
\begin{equation}\label{e54}
    \{x_j + \sigma\sqrt{2\Delta t}\, e_{m}, m = 1,\ldots, M, x_j \in \mathcal{S}\cap \mathcal{D}\} \subset \mathcal{D},
\end{equation}
which is easy to achieve since the diffusivity $\sigma$ is a bounded.

\subsection{The fully discrete scheme}\label{sec:full}
The fully discrete scheme is defined by neglecting all truncation errors and by performing an iterative update from $t_{N_t}$ to $t_0$.

\begin{scheme}\label{s4:full}
Given the temporal spatial partition $\mathcal{T}\times \mathcal{S}$, the terminal condition $u^{N_t}(x_j)$ for $x_j \in \mathcal{S}$. For $n = N_t-1, \ldots, 0$, the approximation of $u(t_n, x)$ is constructed via the following steps:
\begin{itemize}\itemsep0.0cm
\item Step 1: Generate quadrature abscissae in Eqs.~\eqref{E1}, \eqref{E2}, for $x_j \in \mathcal{S} \times \mathcal{D}$
\[
 x_j  + \sigma \sqrt{2\Delta t}e_m, \text{ for } m = 1, \ldots, M, \quad \text{and} \quad x_j + c(a_l),\text{ for } l = 1, \ldots, L.
\]
\item Step 2: Evaluate $u^{n+1,p}(x)$, defined in Eq.~\eqref{eq_interp_E}, at the quadrature abscissae.
  \item Step 3: Compute the approximate expectations 
    \[
    \widehat{\mathbb{E}}\left[u^{n+1,p}(\widehat{X}^{t_n}_{t_{n+1}})\, \big|\, N_{\Delta t} = 0\right] \;\; \text{ and }\;\;  \widetilde{\mathbb{E}}\left[u^{n+1,p}({V}^{t_n}_{t_{n+1}})\, \big|\, N_{\Delta t} = 1\right]
    \]
    via the quadrature rules in Eq.~\eqref{E1} and Eq.~\eqref{E2}, respectively.
    
    \item Step 4: Compute the nodal values $u^{n+1}_j$ through
    \[
    u^{n}_j =\mathbb{P}(N_{\Delta t} = 0) \widehat{\mathbb{E}}\left[u^{n+1,p}(\widehat{X}^{t_n}_{t_{n+}})\, \big|\, N_{\Delta t} = 0\right] +\mathbb{P}(N_{\Delta t} = 1) \widetilde{\mathbb{E}} \left[  u^{n+1,p}({V}_{t_{n+1}}^{t_n})\, \big|\,N_{\Delta t} = 1\right]
    \]
    that is the approximation of the nodal values $u(t_{n+1}, x_j)$ for the interior grid points $x_j \in \mathcal{S} \cap\mathcal{D}$, where  $\widehat{X}_{t_{n+1}}^{t_n}$ and ${V}^{t_n}_{t_{n+1}}$ start from $(t_{n},x_j)$.
\item Step 5: construct the interpolant $u^{n,p}(x)$ by substituting $\{u^n_j\}_{j=1}^J$ into Eq.~\eqref{eq_interp_E}.
\end{itemize}
\end{scheme}

\section{Numerical tests}\label{sec:ex}
In this section, we present a couple of numerical examples to investigate the exit time probability of the L\'evy process $X_t^\alpha$ defined by
\begin{equation}
X_t^\alpha = X_0^\alpha + \int_{0}^{t} (D_{\alpha})^{1/\alpha} dL_s^{\alpha}, \quad \text{with} \quad X_0^\alpha \in \mathcal{D} \subset \mathbb{R}^d.
\end{equation}
Our focus is on solving the exit problem within a finite spatial domain $x \in [0, L]^d$ over a time interval $t \in [0, T]$.

To simplify the analysis and reduce the number of parameters, we introduce a dimensionless formulation by rescaling the spatial and temporal variables
\begin{equation}
\widetilde{x} = {x}/{L}, \quad \widetilde{t} = {t}/{T}.
\end{equation}
Under this rescaling, the problem reduces to computing the exit time probability of the dimensionless L\'evy process $\widetilde{X}_t^\alpha$, defined as
\begin{equation}\label{eq:chi}
\widetilde{X}_t^\alpha = \widetilde{X}_0^\alpha + \int_{0}^{t} \chi  dL_s^{\alpha},
\end{equation}
where $\chi = \frac{(D_{\alpha}T)^{1/\alpha}}{L}$ characterizes the rescaled jump intensity, and the domain is normalized to $(t, x) \in \mathcal{T} \times \mathcal{D} = [0, 1] \times [0, 1]^d$. This transformation balances spatial and temporal scales with the diffusivity parameter, allowing the problem to be described using a single effective diffusivity parameter.
For notational simplicity, we omit the tilde symbol $\widetilde{\cdot}$ in the remainder of the numerical examples.

In the numerical examples, the model is governed by two key physical parameters: the stability index $\alpha \in [1, 2]$ and the dimensionless diffusivity parameter $\chi$. The parameter $\alpha$ determines the degree of nonlocality in the process, with smaller values indicating heavier tails and more distinct jumps. The parameter $\chi$ is defined as $\chi = L_\alpha / L$, where $L_\alpha = (D_\alpha T)^{1/\alpha}$ represents the characteristic spatial scale of fractional diffusion over the time interval $T$. Thus, $\chi$ quantifies the extent of fractional diffusion relative to the system size $L$ during the time scale $T$.
For example, when $\chi \gg 1$, it corresponds to a regime where the domain is small $(L \ll 1)$, the fractional diffusivity is large $(D_\alpha \gg 1)$, and/or the integration time is long $(T \gg 1)$ -- all of which contribute to a dominant nonlocal transport effect across the domain.

\subsection{A 1D symmetric $\alpha$-stable process}\label{sec:ex1}
We consider the dynamics are governed by a one-dimensional symmetric $\alpha$-stable ($S \alpha S$) process $dX_t^{\alpha} = \chi dL_s^\alpha$. Here, $\chi$ is the \textit{normalized jump scale}, previously defined in Eq.~\eqref{eq:chi}, which characterizes the relative strength of fractional diffusion.  We focus on the exit probabilities $P_X(t,x)$ defined in Eq.~\eqref{e2} within the bounded domain $(t,x) \in [0,1] \times [0,1]$.

As the analytical solution for this problem is unknown, we use the direct Monte Carlo (DMC) method with sufficient samples to serve as the ground truth. Specifically, we employ $N_{\rm sample} = 20,000$ samples for the DMC method, ensuring that the reference solution is fully converged, meaning it remains unchanged even with an increase in the number of samples.

\begin{figure}[h!]
    \centering
   {\includegraphics[scale = 0.35]{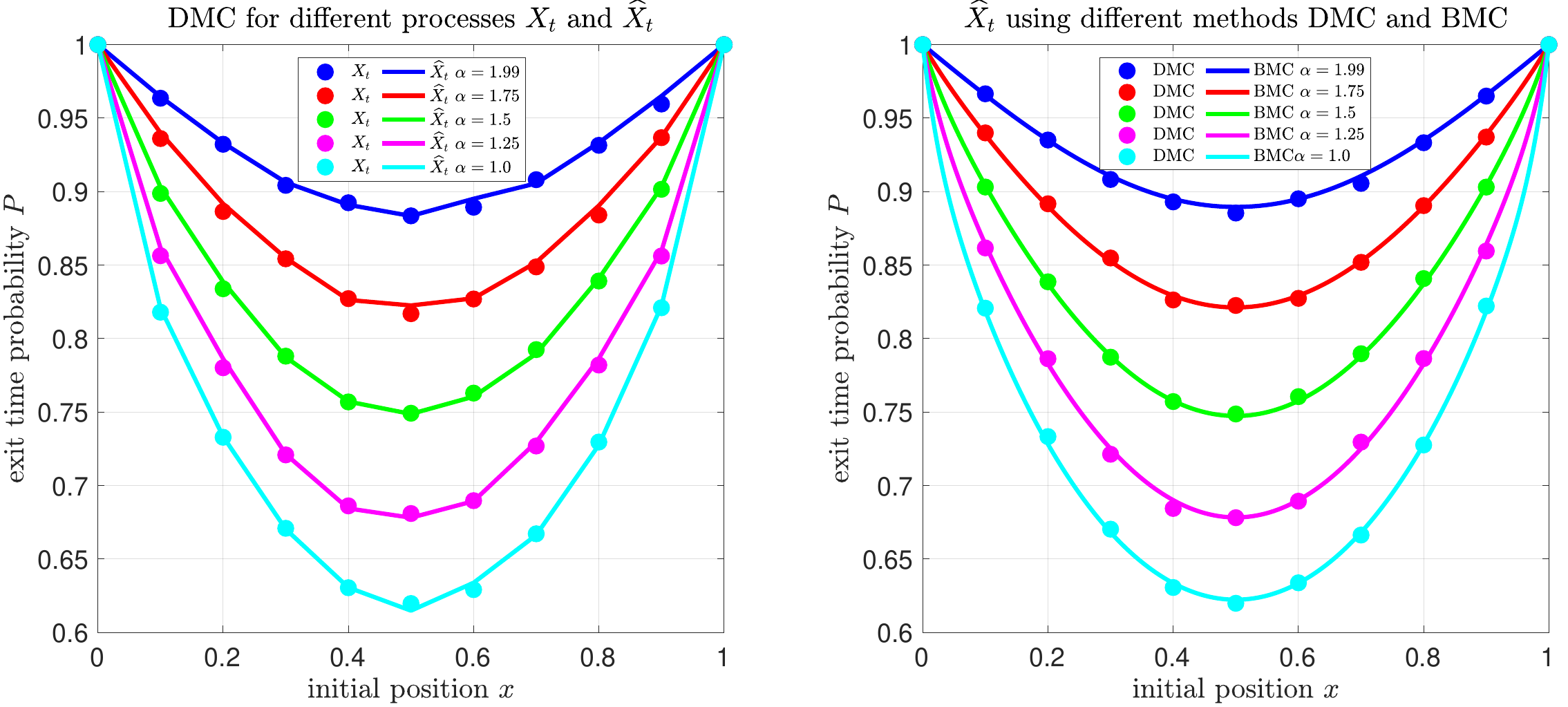} }
    \caption{Left panel: A comparison between $P_X(t,x)$ and $P_{\widehat{X}}(t,x)$ at $t=1$, with both values are computed using the direct Monte Carlo method (DMC). This plot demonstrates $P_{\widehat{X}}(t,x)$ serves as a reliable approximation for $P_X(t,x)$. Right panel: A comparison between DMC (ground truth) and BMC (our method) in solving $P_{\widehat{X}}(t,x)$. A good agreement between our method and ground truth.
    }
    \label{PRE_DMC_BMC}
\end{figure}

We denote the approximation of $X_t^{\alpha}$ in Eq.~\eqref{eq_xtilde} as $\widehat{X}^\alpha_t$. First, we demonstrate that the escape time probability $P_{\widehat{X}}(t,x)$ serves as a reliable approximation for $P_X(t,x)$. 
Utilizing the direct Monte Carlo method (DMC), we compute both $P_X(t,x)$ and $P_{\widehat{X}}(t,x)$, as depicted in the left panel of Fig.~\ref{PRE_DMC_BMC}. It displays the exit time probability $P(t,x)$ for initial positions $x\in [0,1]$ at $t=1$, for $\chi = 0.5$ and $\alpha = 1.99$, $1.75$, $1.5$, $1.25$ and $1.0$. Notably, a remarkable agreement is observed between $P_X(t,x)$ and $P_{\widehat{X}}(t,x)$. The exit time probability tends to increase as the nonlocal parameter $\alpha$ increases since $\chi = 0.5$ refers to a relatively small spatial domain $L$. In which case, the local diffusion tends to higher exit time probability compared to nonlocal diffusion. 
Following this, we move on to a comparative analysis between the direct Monte Carlo method (DMC) and our proposed numerical approach (we denote it as `BMC' since the numerical scheme \ref{s4:full} is a backward algorithm from $t = T_{\rm max}$ to $t=0$).  The right panel of Fig.~\ref{PRE_DMC_BMC} shows the calculation of $P_X(t,x)$, revealing a notable agreement between the outcomes derived from DMC and BMC methods. This observed consistency serves as numerical evidence, reaffirming the precision of our approach.

\begin{figure}[h!]
    \centering
   {\includegraphics[scale = 0.35]{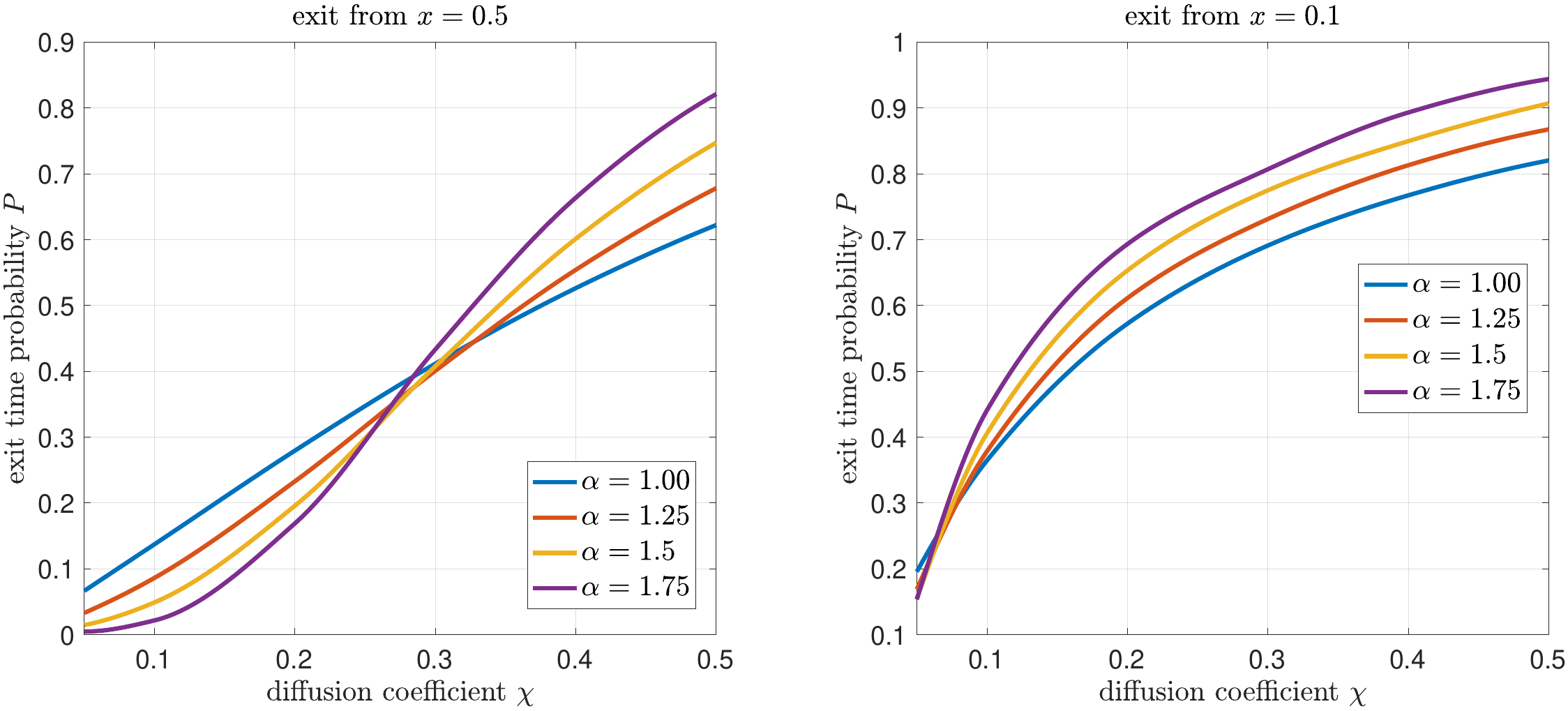} }
    \caption{The exit time probability $P_{\widehat{X}}(t,x)$ at $t=1$ as a function of the diffusion coefficient $\chi$. The left panel corresponds to the initial position $x=0.5$, and the right panel to $x=0.1$. We observe a shift from local to nonlocal diffusion behavior as $\chi$ decreases. Additionally, near-boundary initial positions exhibit faster exits under local diffusion conditions.}
    \label{PRE_chi}
\end{figure}

Next, we examine the influence of the diffusion coefficient $\chi$ on transport behavior by computing the exit time probability $P_{\widehat{X}}(t, x)$ at $t = 1$ across a range of $\chi \in [0.05, 0.5]$. The results are presented in Fig.~\ref{PRE_chi}, showing outcomes from two different initial positions: $x = 0.5$ (left panel) and $x = 0.1$ (right panel).
In the left panel, which corresponds to particles initially located at the center of the domain, we observe a clear trend: as $\chi$ decreases—implying an increase in domain size $L$ or a decrease in the effective fractional diffusion scale $(D_\alpha T)^{1/\alpha}$—the dominant transport behavior transitions from local to nonlocal diffusion. Among the four values of $\alpha$ considered ($\alpha = 1.0$, $1.25$, $1.5$, and $1.75$), the process with $\alpha = 1.75$ shows the fastest exit times for larger $\chi$ values and the slowest for smaller $\chi$ values. This behavior aligns with theoretical expectations: local diffusion is efficient within confined regions, while nonlocal diffusion enables long-range jumps that distribute particles more broadly across the domain.
Furthermore, the exit time probability is sensitive to the initial position of the particle. When the initial position is close to the boundary (as in the right panel, $x = 0.1$), local diffusion tends to result in faster exits—unless $\chi$ is sufficiently small, in which case nonlocal effects dominate.

\begin{figure}[h!]
    \centering
   {\includegraphics[scale = 0.5]{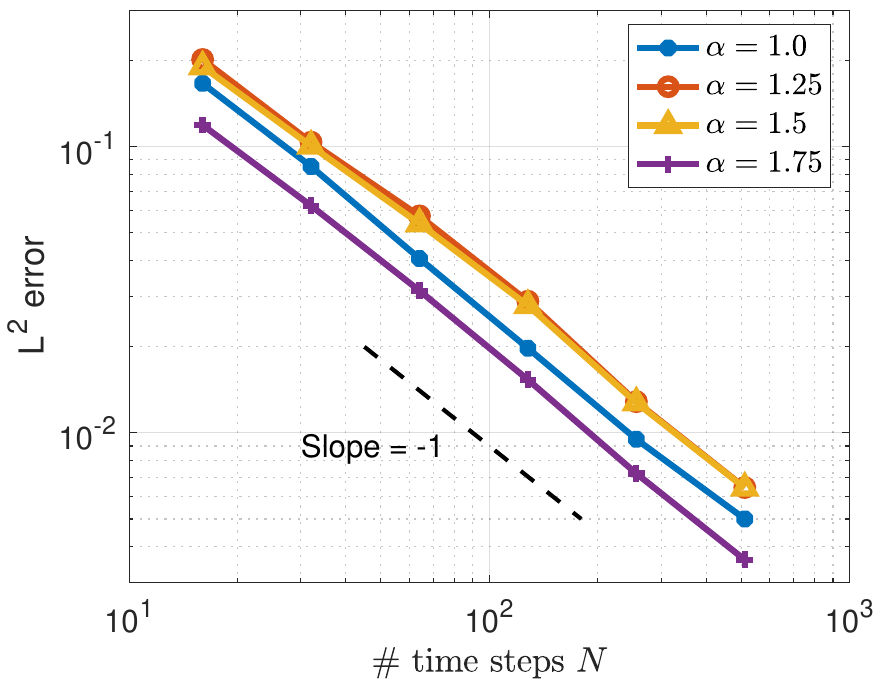} }
    \caption{Error decay in the computation of $P_{\widehat{X}}(t,x)$ at $t=1$ for $\alpha = 1.0$, $1.25$, $1.5$ and $1.75$. The error is measured by the $l^2$-norm of the difference between BMC result and the ground truth (computed by DMC). The proposed method achieves the first-order convergence rate in $\Delta t$.}
    \label{ex1_convergence}
\end{figure}

Finally, we numerically demonstrate the first-order temporal convergence of our method. Fig.~\ref{ex1_convergence} shows the decay of error in the computed exit time probability as the number of time steps $N_t$ increases (i.e., as $\Delta t$ decreases). The error is measured using the $L^2$-norm between the numerical solution obtained via the BMC method and the ground truth  (DMC for $P_{\widehat{X}}(t,x)$). 
As illustrated in Fig.~\ref{ex1_convergence}, the numerical error consistently decreases with smaller time steps, confirming that the method achieves first-order convergence across various values of $\alpha$. For the simulations presented in Figs.~\ref{PRE_DMC_BMC} and \ref{PRE_chi}, we adopt a fine time step of $\Delta t = 10^{-4}$ (corresponding to $N_t = 10^4$) to ensure high accuracy and minimize numerical artifacts.

\subsection{A two-dimensional anisotropic problem}

In this example, we consider a two-dimensional anisotropic transport problem. The transport in the radial direction $ r $ is governed by nonlocal dynamics modeled by a symmetric $\alpha$-stable Lévy process $ L_t^{\alpha} $, while the angular transport in $ \theta $ is locally driven by a standard Brownian motion $ W_t $.
To study the interplay between local and nonlocal anisotropy, we define the following stochastic differential equation on the domain $ (\theta, r) \in \mathcal{D}:= [-\pi, \pi] \times [0, 1] $:
\begin{equation}\label{ex2_sde}
\left\{
\begin{aligned}
d\theta &= P_e\, v_{\theta} (\theta,r) \, dt + \chi\, dW_t,\\ 
dr &= P_e\, v_{r} (\theta,r) \, dt + \chi \, dL_t^{\alpha},
\end{aligned}
\right.
\end{equation}
Here, $\chi$ is the \textit{normalized diffusion scale}, previously defined in Eq.~\eqref{eq:chi}, which characterizes the relative strength of diffusion. In this example, we set $\chi = 0.1$, corresponding to a regime with either a large spatial domain or relatively weak fractional diffusivity.

The velocity field $\mathbf{v} = v_{\theta} \mathbf{e}_{\theta} + v_r \mathbf{e}_r$ is incompressible and defined as:
\[
v_{\theta} = -m\pi \cos(m\pi r)\sin(n\theta), \quad v_r = n \sin(m\pi r)\cos(n\theta).
\]
The model involves three main parameters: the P\'eclet number $P_e$, and the mode numbers $m$ and $n$. The regime $P_e \gg 1$ ($P_e \ll 1$) corresponds to an advection-dominated (diffusion-dominated) regime. In this work, we investigate the impact of $P_e$ on the exit time probability by setting $P_e = 0.1, 10$. The parameters $m$ and $n$ determine the number of convective cells: $m$ specifies the number of cells in the radial direction, while $2n$ controls the number of cells in the angular direction. In this work, we fix $n = 2$ and vary $m = 1, 2$.

\begin{figure}[h!]
    \centering
   {\includegraphics[width=0.9\linewidth]{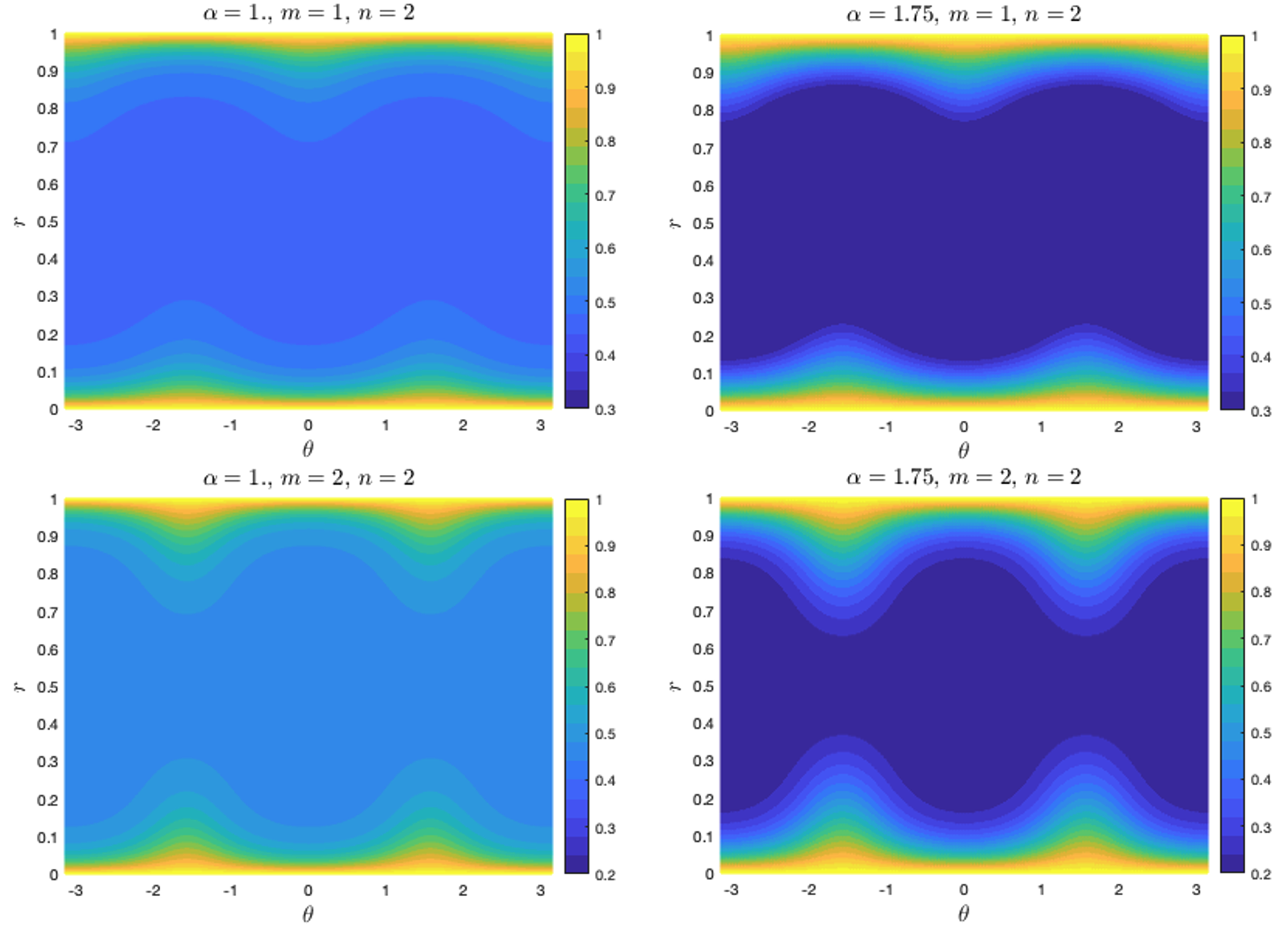} }
 \caption{Exit time probability $P_{\widehat{X}}(t, \theta, r)$ at $t = 1$ for $P_e = 0.1$, and different values of $\alpha=1, 1.75$. Top row: $m = 1$, $n = 2$ with $\alpha = 1$ (left) and $\alpha = 1.75$ (right). Bottom row: $m = 2$, $n = 2$ with $\alpha = 1$ (left) and $\alpha = 1.75$ (right). In the diffusion-dominated regime, particles farther from boundaries exhibit lower exit probabilities, with stronger nonlocality ($\alpha = 1$) leading to higher exit probabilities in the central domain region.}
    \label{ex2fig}
\end{figure}
\begin{figure}[h!]
    \centering
   {\includegraphics[width=0.5\linewidth]{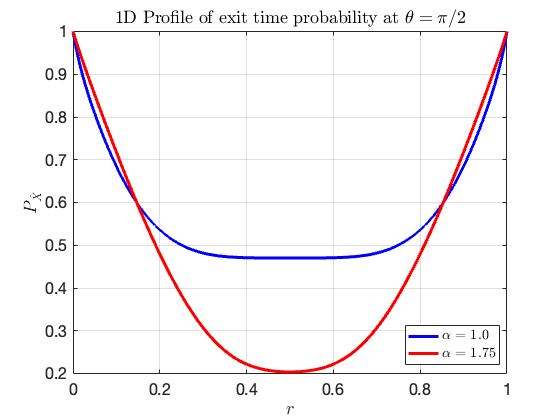} }
\caption{One-dimensional profile of exit time probability $P_{\widehat{X}}(t, \theta, r)$ at $t = 1$ and $\theta = \pi/2$ for $P_e = 0.1$, $m = 2$ and $n = 2$. The comparison between $\alpha = 1$ (blue line) and $\alpha = 1.75$ (red line) demonstrates the transition from local to nonlocal diffusion dominance: near boundaries ($r \approx 0, 1$), local diffusion is more effective ($\alpha = 1.75$ shows higher exit probability), while in the central region, nonlocal effects dominate ($\alpha = 1$ shows higher exit probability).}
    \label{ex2fig_2_1d}
\end{figure}

The exit time probability represents the likelihood that a passive tracer particle starting at a point $(\theta,r)$ reaches either boundary $r = 0$ or $r = 1$ within a given time. The boundary conditions are periodic in $\theta$ and absorbing in $r$, i.e., $P_X(t, \theta, 0) = P_X(t, \theta, 1) = 1$.
Fig.~\ref{ex2fig} presents the exit time probability $P_X(t, \theta, r)$ at final time $T_{\max}=1$ for $P_e = 0.1$, where diffusion dominates the stochastic dynamics. Due to the small $P_e$ value, the weak advection does not generate pronounced convection cells. Consequently, particles located farther from the boundaries exhibit lower exit probabilities. As shown in the figure, when $\chi = 0.1$, which corresponds to a relatively large spatial domain $\mathcal{D}$, the strong nonlocality of the L\'evy process with $\alpha = 1$ (left panel) leads to higher exit probabilities in the central region of the bounded domain compared to the weakly nonlocal process with $\alpha = 1.75$. 
As additional supporting evidence, Fig.~\ref{ex2fig_2_1d} shows the marginal distribution $P_{\widehat{X}}(t, \theta, r)$ at $t = 1$ and $\theta = \pi/2$ for the case $P_e = 0.1$, $m = 2$, $n = 2$. We observe that when the initial location $r$ is near the boundary, local diffusion becomes more significant (the $\alpha = 1.75$ case exhibits higher exit time probability), while nonlocal effects dominate when $r$ is positioned far from the boundary.

\begin{figure}[h!]
    \centering
   {\includegraphics[width=0.9\linewidth]{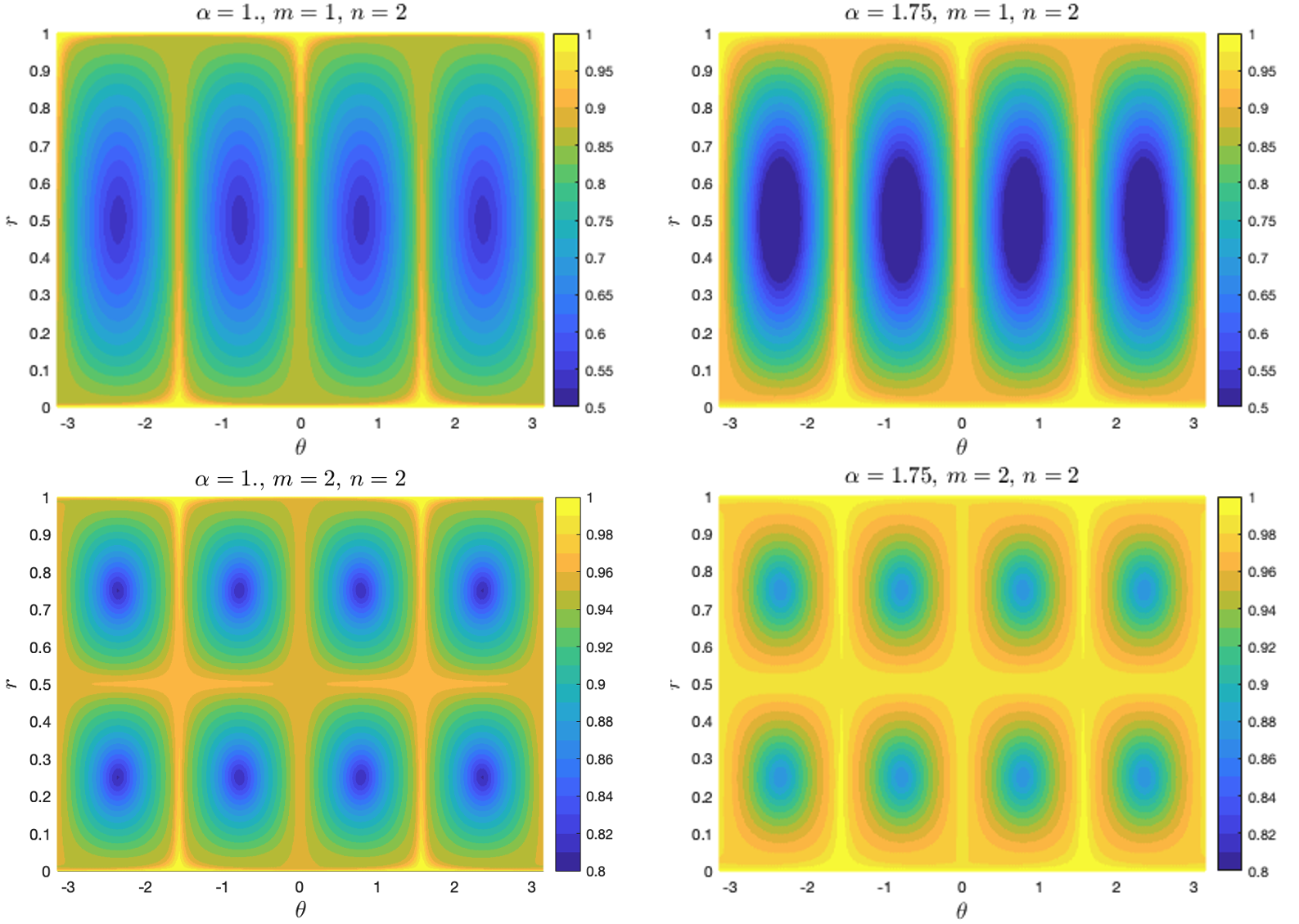} }
  \caption{Exit time probability $P_X(t, \theta, r)$ at $t = 1$ for $P_e = 10$, in the advection-dominated regime. Top row: $m = 1$, $n = 2$ with $\alpha = 1$ (left) and $\alpha = 1.75$ (right). Bottom row: $m = 2$, $n = 2$ with $\alpha = 1$ (left) and $\alpha = 1.75$ (right). The flow exhibits $2mn$ convection cells with elliptic stagnation points (low exit probability) at cell centers and hyperbolic points (high exit probability) along separatrices. For $m = 1$, nonlocal diffusion ($\alpha = 1$) is more effective due to cell centers being far from boundaries. For $m = 2$, local diffusion ($\alpha = 1.75$) dominates as cell centers are closer to boundaries.}
    \label{ex2fig_2}
\end{figure}

When advection effects are considered, i.e., $P_e = 10$, there are $2mn$ convection cells formed. This structure is reflected in the exit time probability: values are lowest near elliptic points (cell centers) and highest near the separatrices and hyperbolic points, indicating faster escape near unstable flow regions. As seen in Fig.~\ref{ex2fig_2}, the flow exhibits clearly defined convection rolls, with elliptic stagnation points at their centers and hyperbolic points along vertical separatrices. When $m = 1$, the four cell centers are located at $(\theta, r) = (\pm\pi/4, 1/2)$ and $(\theta, r) = (\pm 3\pi/4, 1/2)$. When $m = 2$, eight cell centers are positioned at $(\theta, r) = (\pm\pi/4, 1/4)$, $(\theta, r) = (\pm\pi/4, 3/4)$, $(\theta, r) = (\pm 3\pi/4, 1/4)$, and $(\theta, r) = (\pm 3\pi/4, 3/4)$. When $m = 1$, the cell centers are relatively far from the boundary. In this case, nonlocal diffusion is more impactful since it enables long-range jumps that distribute particles more broadly across the domain. From the top row of Fig.~\ref{ex2fig_2}, we observe that the $\alpha = 1$ case exhibits higher exit time probability $P_{\widehat{X}}(t, \theta, r)$ near the cell centers. When $m = 2$, the cell centers are relatively closer to the boundary. In this case, local diffusion is more efficient at moving particles to the boundary. The $\alpha = 1.75$ L\'evy process exhibits higher exit time probability as shown in the bottom row.

\section{Conclusion}
We have developed an efficient probabilistic numerical method to compute the exit time probability associated with stochastic processes driven by symmetric $\alpha$-stable L\'evy noise. By approximating the original $\alpha$-stable process with a combination of Brownian motion and a compound Poisson process, the proposed approach overcomes key computational challenges inherent in nonlocal fractional models. Our method reformulates the exit problem into a terminal-value partial integro-differential equation (PIDE), which is then solved via a probabilistic scheme based on the Feynman-Kac formula and tailored quadrature rules.

The numerical scheme achieves first-order convergence in time and avoids solving large dense linear systems, offering significant efficiency gains over traditional PDE-based and Monte Carlo approaches. Through extensive numerical experiments in one and two dimensions, we demonstrated that the method reliably reproduces key statistical features of the $\alpha$-stable process—such as the heavy-tailed distribution and superdiffusive scaling behavior—and delivers accurate approximations of exit time probabilities across a range of parameters. The two-dimensional anisotropic example further illustrates the method's capability to capture the interplay between local and nonlocal transport mechanisms in complex flow configurations.

Overall, the proposed scheme provides a robust and scalable tool for studying exit phenomena in systems characterized by anomalous diffusion, with potential applications in physics, finance, and other fields involving nonlocal transport processes.

\section*{Acknowledgements}
This material is based upon work supported in part by the U.S.~Department of Energy, Office of Science, Offices of Advanced Scientific Computing Research and Fusion Energy Science, and by the Laboratory Directed Research and Development program at the Oak Ridge National Laboratory, which is operated by UT-Battelle, LLC, for the U.S.~Department of Energy under Contract DE-AC05-00OR22725.

\bibliographystyle{abbrv}
\bibliography{PIDE_ref,library,myref,fractional}

\end{document}